\documentclass[twoside,12pt]{amsart}
\usepackage[left=1.1in,right=1.2in, top=1.2in, bottom=1.5in]{geometry}
\usepackage{amssymb,amsmath,amscd,amsthm,mathrsfs,graphicx} 
\usepackage[numbers]{natbib}
\usepackage{tikz}
\usepackage{blkarray}
\usepackage{xcolor}
\usepackage{upgreek,xspace}
\usepackage[pdftex,breaklinks,colorlinks]{hyperref} 
\hypersetup{linkcolor=blue,citecolor=blue}
\usepackage{cleveref}
\usepackage{caption}
\usepackage[color=red!30!,textwidth=2cm]{todonotes}
\usepackage{blindtext}
\usepackage{datetime}
\newdateformat{monthdayyeardate}{%
  \monthname[\THEMONTH]~\THEDAY, \THEYEAR}
\usepackage{parskip}

\usepackage[normalem]{ulem} 

\usepackage{filemod}

\usepackage{fancyhdr}

\newcommand{\ZZ}{{\mathbb{Z}}}
\newcommand{\CC}{{\mathbb{C}}}

\newcommand{\QQ}{{\mathbb{Q}}}

\renewcommand{\S}{{\mathbb{S}}}

\newcommand{\T}{\mathcal{T}}

\newcommand{\bx}{{\mathbf{x}}}
\newcommand{\by}{{\mathbf{y}}}

\newcommand{\BC}{\mathfrak{C}}
\newcommand{\BT}{\mathfrak{T}}
\newcommand{\BZ}{\mathfrak{Z}}

\newcommand{\GC}{\mathscr{C}}
\newcommand{\GT}{\mathscr{T}}
\newcommand{\GZ}{\mathscr{Z}}

\newcommand{\gb}{\mathfrak{b}}

\newcommand{\gc}{\mathfrak{c}}
\newcommand{\gB}{\mathrm{b}}

\newcommand{\gh}{\mathfrak{h}}

\newcommand{\ot}{\overline{\T}}
\newcommand{\pt}{\partial{\T}}
\newcommand{\dl}{\T_1 \ud \T_2}

\newcommand{\ozz}{\overline{\ZZ }}

 \newcommand{\bus}[1]{\overline{#1}^{\textrm{\tiny Bus}}}
 \newcommand{\bxi}{\pmb{\xi}}

\newcommand{\ud}{\uparrow\downarrow}
 
\renewcommand{\i}{\iota}

 \usepackage{xcolor}

\usepackage{soul}

\setstcolor{blue}

\theoremstyle{plain}
\newtheorem{theorem}{Theorem}[section]

\newtheorem{lemma}[theorem]{Lemma}

\theoremstyle{definition}

\newtheorem{remark}[theorem]{Remark}
\newtheorem{example}{Example}[section]

 \newcommand{\gcom}{(\T_1 \ud \T_2)^{\tiny{{\textrm{HM}}}}}
  \newcommand{\bcom}{(\T_1 \ud \T_2)^{\textrm{Bus}}}

\title{\Large C\MakeLowercase {ompactifications of horospheric products}}
 
\author{\large B\MakeLowercase {ehrang} F\MakeLowercase {orghani and} K\MakeLowercase {eivan} M\MakeLowercase {allahi}-K\MakeLowercase {arai}}

\address[1]{The College of Charleston, 66 George Street, Charleston, SC 29403, USA}
\email{forghanib@cofc.edu}
\address[2]{Jacobs University, Campus Ring 1, 28759 Bremen, Germany}
\email{k.mallahikarai@jacobs-university.de }

\begin{document}

\begin{abstract}
We define and study a new compactification, called {\it the height compactification} of the horospheric product of two infinite trees.  We will provide 
a complete description of this compactification. In particular, we show that this compactification is isomorphic to the Busemann compactification when all the vertices of both trees have degree at least three, which also leads to a precise description of the 
Busemann functions in terms of the points in the geometric compactification of each tree. We will discuss an application to the asymptotic behavior of integrable ergodic cocycles with values in the isometry group of such horospheric product.
\end{abstract}

\keywords{Horospheric product, Busemann compactification}

\maketitle

\section{\bf Introduction}
 
The notion of Busemann compactification of metric spaces, introduced in \cite{Gromov}, is an instrumental tool in studying geometric and probabilistic aspects of groups, graphs, and manifolds. 
This construction appears in different contexts, for instance, the work of Cheeger and Gromoll theorem \cite{Cheeger-Gromoll} in non-negative curvature,  work of Karlsson and Ledrappier \cite{Karlsson-Ledrappier-law} and \cite{Karlsson-Ledrappier-Linear}, and Gou\"{e}zel-Karlsson \cite{Gouzel-Karlsson} in the theory of random walks on groups; see Karlsson \cite{Karlsson} and the references therein for more details. 

A challenging problem in this theory is to identify the points in the Busemann compactification in terms of the  geometry of the underlying space.  Notable examples where such a geometric description is possible include the work of Walsh \cite{Walsh} on the Busemann compactification of the Teichm\"{u}ller space,  Horbez's result \cite{Horbez} on  the Busemann compactification of the Outer space, and the work of Maher-Tiozzo \cite{Maher-Tiozzo} where the Busemann compactification of non-proper hyperbolic spaces has been investigated. 

The goal of this work is twofold. On the one hand, we will provide a general construction of new compactifications for a metric space $X$ starting from an initial compactification $\widetilde{X}$ and a continuous mapping $h: X \to K$, where $K$ is itself a compact metric space. Roughly speaking, this space, called the {\it mapping compactification}, is the closure of the diagonal embedding of $X$ in $\widetilde{X} \times K$, see Section \ref{sec:mapping} for details. 
We will then apply this construction to the specific case of horospheric products of trees. This leads to a compactification that, under appropriate conditions, turns out to be  isomorphic to  the Busemann compactification. Interestingly, these compactifications are both larger than the one obtained by embedding the horospheric product into the product of trees. We will also provide a decomposition of these two compactifications into components, each of which can be precisely described.

Horospheric products of trees were initially brought into the limelight of geometric group theory by Diestel and Leader \cite{Diestel-Leader} as possible candidates for 
vertex transitive graphs that are not quasi-isometric to Cayley graphs. The question of the existence of such graphs were previously posed by  Woess \cite{Wolfgang}.
Eskin, Fisher, and Whyte \cite{Eskin-Fisher-Whyte,Eskin-Fisher-Whyte-1} settled the conjecture by proving that these graphs are indeed not quasi-isometric to Cayley graphs.
 
Ever since their introduction, these graphs have been a source of many challenging questions that involve an interplay between probabilistic and geometrical aspects.
Some of these work include convergence of random walks \cite{Bertacchi},  spectral radius of simple random walk  \cite{Coste-Woess} and \cite{Barthodli-Woess},  Martin boundary and minimal harmonic functions  in  \cite{Brofferio-Woess05} and \cite{Brofferio-Woess06},  the Poisson boundaries of discrete isometry groups of horospheric products \cite{Batholdi-Neuhauser-W} \cite{Kaimanovich-Sobieczky} and  the Poisson boundaries of locally compact isometry groups of horospheric products \cite{Forghani-Tiozzo}. 
It is also worth pointing out that, disparate as they might appear, the Martin compactification, the Busemann 
compactification, and the Thurston compactification of the Teichm\"uller space have the common feature that they are all based on embedding a certain space $X$ 
into a projective space of functions on $X$, and then closing it. For the Martin compactification (Thurston compactification, respectively) this is accomplished by means of the Green kernel (intersection function, respectively). A unified treatment of these compactifications due to Constantinescu-Cornea can be found in \cite{Brelot}.

Let $\T_i$ be locally finite infinite trees where the degree of each vertex is at least $2$ for $i=1,2$. We fix an infinite geodesic $\gamma^i$ and a base point $o^i$ in $\T_i$.  Denote by $\ot_i$ as the geometric compactification of tree $\T_i$, which is isomorphic with its Busemann  compactification  (see Section~\ref{sec:trees} for precise definitions).
Denote by $\gh_i$ the Busemann function (see Section~\ref{sec:trees} for the definition) associated with $\gamma^i$ and base point $o^i$. Denote by  $\dl$ the horospheric product (see Section~\ref{sec:horospheric} for the definition) associated with $(\T_i,o^i,\gamma^i)$ for $i=1,2$. The geometric compactification of $\dl$ is 
defined as the 
closure of $\dl$ as a subset of $\ot_1\times \ot_2$ in the product topology.

Indeed, the geometric compactification of $\dl$ is the mapping compactification under the embedding $\dl \to \ot_1\times \ot_2  \times \{1\}$. As it will be shown, the geometric compactification does not provide any instructive information about level sets of the horospheric products $\dl$. 
We will consider another mapping compactification of $\dl$ via a natural embedding 
$\dl \to \ot_1\times \ot_2 \times \ozz$. More precisely, define the \emph{height compactification} of $\dl$, denoted by 
$\gcom$, 
as the closure of $\dl$, when viewed as a subset of $\ot_1\times \ot_2 \times \ozz$, where the third component records the value of the Busemann functions associated to $\gamma^1$.  We will show that  sequence $\bx_n \in \dl$ converges in $\gcom$ if and only if it converges in the Busemann compactification $\bcom$. This provides an isomorphism of $\gcom$ and $\bcom$, viewed as compactifications of $\dl$.

\begin{theorem}[Main Result]\label{thm:main}
Let $\T_1$ and $\T_2$ be two  locally finite trees where the degree of each vertex is at least three. Then, the Busemann compactification of $\dl$  is isomorphic to the height compactification of
$\dl$. 
\end{theorem}

It is worth mentioning that the boundary of the height compactification is isomorphic to the horocycle boundary of two regular trees introduced by Brofferio and Woess in \cite{Brofferio-Woess05}.
In their construction of horocycle boundary, one first compactifies each tree with a metric that is different from the graph metric, and then the horocycle compactification 
is defined as the closure of embedding of a horospheric product to the product of these compactifications.

Brofferio and Woess  \cite{Brofferio-Woess05} proved that the horocycle boundary of product of two regular trees, when equipped with the hitting measure, can be identified as the Martin boundary of a class random walks. Combining results in \cite{Brofferio-Woess05} with Theorem \ref{thm:main}, we conclude that the Busemann compactification can be identified with the underlying topological space of the Martin boundary of some random walks. The special case where $\T_1$ and $\T_2$ are both $3$-regular (that is, the degrees of all vertices of $\T_1$ and $\T_2$ are equal to $3$), where $\dl$ coincides with a Cayley graph of a lamplighter group has been studied using different methods in \cite{Jones-Kelsey}. One can extend the definition of horospheric products to Gromov hyperbolic spaces. It would be interesting to provide a succinct description of the Busemann compactification of horospheric products of hyperbolic spaces analogous to Theorem \ref{thm:main}. This question has been studied for another compactification of horospheric products of two Gromov hyperbolic spaces in \cite{Ferragut}.

Aside from the isomorphism above, we will also provide a decomposition of both boundaries into subspaces each of which has a simple topological description. A less technical version of these theorems (Theorems~ \ref{thm:geo} and \ref{thm:bus=HM}) can be stated as follows.  

\begin{theorem}\label{thm:decomp-simple}
Let $\T_1$ and $\T_2$ be two   locally finite trees where the degree of each vertex is at least three. Then, the boundary points of the Busemann compactification and the height compactification can each be decomposed as 
$$
C_1 \sqcup C_2 \sqcup T_1 \sqcup T_2 \sqcup Z
 $$
where for $i=1,2$, the space $C_i$ is homeomorphic to the geometric boundary of $\T_i$, the elements of $T_i$ are in bijection with the vertices of $\T_i$ and $Z$ is a discrete sets parametrized by elements of $\ZZ$.  Moreover, the closure of $Z$ intersects each one of $C_1$ and $C_2$ in one point. 
\end{theorem}

Our proof is based on an explicit description of the Busemann functions of the horospheric product of two trees in Section~\ref{sec:Busemann}, which, in turn, relies on the description of the graph metric of $\dl$ due to Bertacchi \cite{Bertacchi} (for regular trees) and  Kaimanovich and Sobieczky \cite{Kaimanovich-Sobieczky} (for all trees). At the end of the paper, we will describe the additional ingredients needed to generalize  Theorem \ref{thm:decomp-simple} to the case of trees which might have vertices of degree $2$; see Theorem \ref{thm:general} and the preceding discussion.

\textcolor{black}{As an application, we will show that Theorem \ref{thm:decomp-simple} combined with  Karlsson-Ledrappier's law of large numbers \cite{Karlsson-Ledrappier-law} provides non-trivial information on asymptotic behavior of $Z_no$ where $Z_n$ is an ergodic 
integrable cocycle with values in the isometry group of $\dl$. See Remark \ref{cocycle} for the precise statement.}

This paper is structured as follows. In Section~\ref{sec:comp-metric } we provide the definition of the notion of compactification for general metric spaces with two examples, mapping and Busemann compactifications. In Section~\ref{sec:trees} we set the notation regarding trees and review a number of basic results regarding their Busemann compactifications. Section~\ref{sec:horospheric} is devoted to definition of horospheric product of two trees and computing its height  compactification.  Section~\ref{sec:Busemann} includes identification of Busemann compactification and the main result.

\subsection*{Acknowledgements} 
We would like to thank V. Kaimanovich for detailed comments on an initial version of this paper and bringing the connection between various compactifications through work of Constantinescu-Cornea with reference \cite{Brelot} to our attention. 
We would also thank G. Tiozzo and W. Woess for pointing to results and references including identification of Busemann compactifications.
Finally, we would like to thank the referees for careful reading of the manuscript and detailed remarks.

\section{\bf Compactifications of metric spaces} \label{sec:comp-metric }
Let $(X,d)$ be an unbounded metric space. A sequence $(x_n)$ in $X$ is called divergent if for some (equivalently, any) base point $o \in X$, we have $d(x_n, o) \to \infty$.  A sequence $(x_n)$ in $X$ is called {\it eventually constant} if there exists $k \ge 1$ such that $x_n=x_m$ for all $m,n\geq k$. The value of $x_k$ will be called the {\it eventual value} of the sequence. A convergent sequence in a discrete set is easily seen to be eventually constant.

By a compactification of $X$ we mean a pair $( \widetilde{ X}, \i)$ 
where $ \widetilde{ X}$ is a compact metric space and $ \i: X \to \widetilde{ X}$ is an injection such that $\i( X)$ is dense in $ \widetilde{ X}$. An infinite discrete space has many compactifications. The one-point compactification and the Stone-\v{C}ech compactification are two prominent examples. Compactifications $( \widetilde{ X}_1, \i_1)$ and $( \widetilde{ X}_2, \i_2)$ are
said to be isomorphic if there exists a homeomorphism $ \phi : \widetilde{ X}_1 \to \widetilde{ X}_2$ such that $ \phi \circ \i_1= \i_2$. This is equivalent to the condition that  for a sequence $(x_n)$ in $X$ we have $\i_1(x_n) \to \xi_1 \in \widetilde{ X}_1$ iff 
$\i_2(x_n) \to \phi(\xi_1) \in \widetilde{ X}_2$.

\subsection{Mapping compactification}\label{sec:mapping} Suppose $X$ is a metric space equipped with a continuous map $h: X \to K$, where $K$ is a compact metric space. Suppose  $( \widetilde{ X}, \i)$  is some compactification of $X$. Consider the map
$$ \widehat{h}: X \to \widetilde{ X} \times K, \quad  \widehat{h}(x)= (\i(x), h(x)).$$
The closure of $ \widehat{h}(X)$ in $ \widetilde{ X} \times  K$ is called the {\it mapping compactification} of $X$ with respect to $h:X \to K$ and compactification $\widetilde{ X} $. 
We will now give a few examples of mapping compactification. 
\begin{example}
Let $K$ be one point, then mapping compactification of $X$ is isomorphic to $ \widetilde{ X}$. 
\begin{example}
Consider $\Bbb Z$ with one point compactification. Set $K=\{-1,1\}$ and let $h$ be the sign function on integers. Then the corresponding mapping compactification is the two point compactification on $\Bbb Z$. 
\end{example}

\begin{example}
Consider $\Bbb R$ with one point compactification. Set $K=[0,1]$ and $h(x)=\frac{e^x}{1+e^x}$. Then the mapping compactification  with respect to $h$ would be the two points compactification of $\Bbb R$.
\end{example}

\end{example}

 \begin{example} \vspace{3mm}
Let $\ZZ$ be the set of integers equipped with the discrete topology. Set $$\S^1= \{ z \in \CC: |z|=1 \}$$ to denote the unit circle in the complex plane. For an irrational $ \alpha \in [0,1]$, define the map $e_{ \alpha}: \ZZ \to \S^1$ by $ e_{ \alpha} (n)= \exp (2 \pi i \alpha n )$. It is well known that $ e_{ \alpha}$ 
is injective and that $ \overline{e_{ \alpha}( \ZZ)}= \S^1$. Viewed in this way, $\S^1$ can be regarded as a compactification of $\ZZ$; denote this compactification by $\overline{\ZZ}_{ \alpha}$.  Now suppose that $\beta$ is a real number such that $ 1, \alpha, \beta $ are linearly independent over $\QQ$. Consider the map
\[ h_{ \alpha, \beta}: \ZZ \to \overline{\ZZ}_{ \alpha} \times \S^1, \qquad h_{ \alpha, \beta}( n)= (e_{ \alpha}(n), e_{\beta}(n) ). \]
It follows from Kronecker's theorem that the image of $h_{ \alpha, \beta}$ is dense in $\S^1 \times \S^1$. As a result, we can regard
$ \S^1 \times \S^1$ as a compactification of $\ZZ$. On the other hand, when 
$ \{ 1, \alpha, \beta \}$ are linearly dependent, the closure of $h_{ \alpha, \beta}(\ZZ)$ is still homeomorphic to $\S^1$. 
As a corollary, we also see that when $ \{ 1, \alpha, \beta \}$ is linearly independent, then the compactifications $\overline{\ZZ}_{ \alpha}$ and $\overline{\ZZ}_{ \beta}$ are not isomorphic. 
Similar constructions can be made using other compact topological groups that contain 
elements of infinite order. If $G$ is such a group, and $g \in G$ is an element of infinite order in $G$, then the closure of the cyclic group generated by $g$ is a compactification of $\ZZ$. Such groups are called \emph{procyclic groups}. 
 \end{example}

Other examples will be provided later; see Theorem~ \ref{thm:geo}.

\subsection{The Busemann compactification}
Fix an arbitrary base point $o \in X$. Denote by $C(X)$ be the space of all real-valued continuous functions on $X$, topologized by the uniform convergence on compact sets topology. 
We will define a map $\Phi:X \to C(X)$ as follows. For $z \in X$, set 
$$
z \mapsto \gb_z(x) := d(z,x) - d(z,o).
$$
Note that $|\gb_z(x) |\leq d(x,o)$, implying that $ \gb_z$ is a continuous function, and for a fixed $x \in X$, the set $\{\gb_z(x)  \}_{ x \in X}$ has a uniform upper bound.  
Note that each 
$ \gb_z$ is a $1$-Lipschitz function with respect to the metric $d$. These allow one to embed 
 $X$ into the compact space $\prod_{x \in X} [- d(x,o), d(x,o)]$. 
Tychonoff's theorem guarantees that $\overline{\Phi(X)}$ is compact. We will call $\overline{\Phi(X)}$ the Busemann compactification of $X$ and denote it  by 
$\bus{X}$. This compactification was introduced by Gromov in \cite{Gromov}. Let us remark in passing that the topology defined in \cite{Gromov} can sometimes be different from the one considered here. See \cite[p. 7]{Karlsson} for a discussion on this subtle point.

It is sometimes useful to work with the Busemann cocycle. Fixing 
point $y, z \in X$, it is easy to see that if $z_n$ converges to $\gamma$ in the Busemann compactification of $X$, then the limit
\[ z_n \mapsto d(y, z_n) - d(x,z_n) \]
also exists. This limit is denoted by $\upbeta_{\gamma}(x,y)$. It is easy to see that $\upbeta_\gamma$ satisfies the cocycle identity:
\[ \upbeta_\gamma(x, y)+ \upbeta_\gamma(y, z)= \upbeta_\gamma(x, z). \]
 
\section{\bf Compactifications of trees}\label{sec:trees} In this section we will define two different compactifications of infinite trees and prove a number of elementary but useful properties for them. 
All tress considered in this section are assumed to be  locally finite, that is, they have the property that the degree of each vertex is finite.  

Let $\T$ be a tree with a fixed base point $o$. We will always assume that the degree of every vertex of $\T$ is at least $2$. An infinite geodesic in $\T$ is  a non-backtracking infinite path emanating from $o$; more precisely, an infinite geodesic is an infinite sequence $\xi=(z_n)_{n\geq0}$ in $\T$ such that $z_0=o$, and for all $i \ge 0$
$z_i$ and $z_{i+1}$ form an edge in $\T$ and  the set $\{ z_i, z_{i+1}, z_{i+2} \}$ consists of three distinct points. 
Define the geometric boundary of $\T$ as the space of all infinite geodesics in $\T$. We denote the geometric boundary of $\T$ by $\pt$.  
By $[y,x]=[x_0=y,x_1,\hdots,x_{n-1},x_n=x]$ denote the unique path between vertices $y$ and $x$ in $\T$.  Consider finite paths $[o,x]=[x_0,\hdots,x_n]$ and $[o,y]=[y_0,\hdots,y_n]$ in the tree $\T$. Let $N(x,y)$ be the length of the longest common path  starting from $o$ in the paths $[o,x]$ and $[o,y]$, more precisely, 
$N(x, y)= \max \{  0 \le k\le n : x_k = y_k \}$.  For an infinite geodesic $\xi=(z_n)_{n \ge 0}$ and the vertex $x$ in the tree $\T$, let $N(x,\xi) : = N (x, z_k)$ where $k= d(o,x)$. Denote by  $[o,x]_\xi= z_{N(x,\xi)}$  the unique common vertex in the finite path $[o,x]$ and the infinite geodesic $\xi$ with the longest distance from the base point $o$.   We say a sequence of $(x_n)$ in $\T$ converges to the boundary point $\xi=(z_n)$ in $\pt$ if and only if  $N(x_n, \xi)=d(o,[o,x_n]_\xi) $ goes to infinity as $n$ approaches infinity.
Similarly, for two distinct infinite geodesic $\xi_1$ and $\xi_2$, we can define the $N(\xi_1,\xi_2)$ as the distance of the longest common path in the geodesic $\xi_1$ and $\xi_2$.  The boundary of the tree can be equipped with the following metric  $$d( \xi_1, \xi_2)= 2^{-N(\xi_1, \xi_2)}.$$ 
With defined topology on $\T$ and $\pt$, we can see $\T\sqcup\pt$ is
a compactification of $\T$ with respect to the graph metric on $\T$.  We call $\ot =\T\sqcup\pt$  the geometric compactification of the tree $\T$. It is worth mentioning that the geometric compactifications of a tree is  isomorphic to its Gromov compactification,   viewed as a hyperbolic space. 

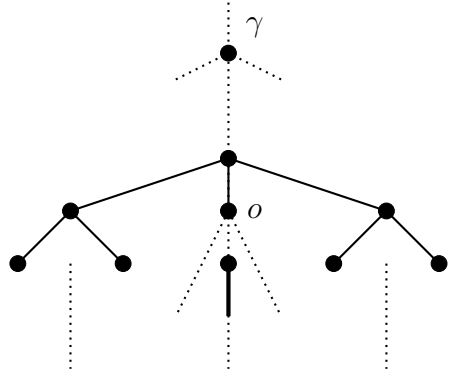
\begin{figure}[h!]
\begin{center}
\begin{tikzpicture}[scale=0.7]
         \draw[fill,black] (0,2) circle(.15cm);
        \foreach \x in {-3,0,3}
        \draw[fill] (\x,1) circle(.15cm);
            \foreach \x in {-3,0,3}
        \draw[black, thick] (\x,1)--(0,2);
        \foreach \x in {-4,-2,0,2,4}
        \draw[fill] (\x,0) circle(.15cm);
        \draw[fill,ultra thick] (0,0)--(0,-1);
        \foreach \x in {-1,1}
        \draw[dotted, thick] (\x,3.5)--(0,4);
        \draw[dotted, thick] (0,0)--(0,5);
         \draw[fill] (0,4) circle(.15cm);
          \foreach \x in {-3,0,3}
            \draw[dotted, thick] (\x,0)--(\x,-2);
            \node at (0.5,1) {$o$};
             \node at (0.5,4.5) {$\gamma$};
            
             \foreach \x in {-1,1}
            \draw[dotted,thick] (0,1)--(\x,-1);
        
            \foreach \x in {-4,-2}
        \draw[fill,  thick] (\x,0)--(-3,1);
              \foreach \x in {4,2}
        \draw[fill,thick] (\x,0)--(3,1);
\end{tikzpicture}
\end{center}
\caption{An example of a rooted tree with base point $o$ and a boundary point $\gamma$.}
\end{figure}

Let $\T$ be a tree with a base point $o$, and let $\gamma=(z_n) \in \pt$.  
Set 
$$ 
\gb_\gamma(y) = d ( y, [o,y]_\gamma) - d(o,[o,y]_\gamma).
$$

We will call $ \gb_\gamma$ the height function associated to $\gamma$. 
For $k \in \ZZ$, we will denote by $H_k(\T)$ the level set of $\gb_\gamma$ defined by $H_k =\{ y \in \T\ :\ \gb_\gamma(y) = k\}$. 
For $k \in \ZZ$, define $\gc_k: \T \to \ZZ$ by 
$$
\gc_k(x) = |k|-| k-\gb_\gamma(x) |.
$$

It is known that the geometric compactification of   $\T$ is isomorphic to the Busemann compactification of $\T$ with respect to its graph metric. We will refer to this fact later when proving the main result.  
\begin{lemma}\label{lem:geometric-busemann-tree}
The Busemann compactification and the geometric compactification of a tree are isomorphic. Moreover, a sequence $(x_n)$ in $\T$ convergence to the boundary point $\gamma$ in the geometric boundary if and only if  for every $y$ in $\T$
$$
\lim_n \gb_{x_n} (y) = \gb_\gamma(y)
$$
The Busemann cocycles can be formulated based on the height function:
$$
\upbeta_{\gamma}(x,y) = \gb_\gamma(y) - \gb_\gamma(x).
$$
\end{lemma}

The following two lemmas will be used later.

\begin{lemma}\label{lem:bounded-tree}
Let $x_n \in \T$ be a divergent sequence of points with the property that $ \gb_\gamma(x_n)$ is bounded. Then $x_n \to \gamma$
in the geometric boundary of $\T$.  
\end{lemma}
\begin{proof}
We have 
\[ |d( o, [o,x_n]_\gamma )- d( x_n, [o,x_n]_\gamma )|= |\upbeta_{\gamma}(o, x_n)| =|\gb_\gamma(x_n)| \le C. \]
Also note that 
\[ d( o, [o,x_n]_\gamma )+  d( x_n, [o,x_n]_\gamma ) \ge d( o, x_n) \to \infty. \]
This implies that $ d( o, [o,x_n]_\gamma ) \to \infty$, which is equivalent to  $x_n \to \gamma$ in the geometric boundary. 
\end{proof}
\begin{lemma}\label{lem:constant-cocycle}
Let  $(x_n)$ be a sequence  in $\T$ with the property that $ \gb_\gamma(x_n)$ is unbounded. If $x_n$ converges to a boundary point in $\pt$, then for every vertices $x$ and $y$ in $\T$, 

$$
\lim_n \Big[|\upbeta_\gamma(x,x_n)| - | \upbeta_\gamma(y,x_n) | \Big] \in \Big\{\upbeta_\gamma(x,y),  \upbeta_\gamma(y,x)\Big\}.
$$
\end{lemma}
\begin{proof}
Because $x_n$ converges to a boundary point, there are two cases:
\begin{enumerate}
\item$ The case  \upbeta_\gamma(o,x_n) =\gb_\gamma(x_n) \to - \infty$: This is equivalent to $x_n \to \gamma$. In this case, for  sufficiently large values of $n$ the values 
$\upbeta_\gamma(x,x_n) = \upbeta_\gamma(o,x_n) - \upbeta_\gamma(o,x)$ and  $\upbeta_\gamma(y,x_n) = \upbeta_\gamma(o,x_n) - \upbeta_\gamma(o,y)$ are both negative. This implies that 
$$
\lim_n \Big[ |\upbeta_\gamma (x_n,x)| - | \upbeta_\gamma (x_n,y) | \Big]=  \upbeta_\gamma(y,x)=\gb_\gamma(x)-\gb_\gamma(y).
$$

\item The case $ \upbeta_\gamma(o,x_n)= \gb_\gamma(x_n)  \to + \infty $ is dealt with in a similar fashion. The only difference is that the same quantities will be eventually positive as $n \to \infty$ and one obtains $\upbeta(x,y)=\gb_\gamma(y) - \gb_\gamma(x)$ as the limit. 

\end{enumerate}

\end{proof}

\section{\bf Horospheric products of trees}\label{sec:horospheric}
 In this section we will recall the construction of horospheric product of rooted trees.  Let $\T_1$ and $\T_2$ be two  rooted trees with base points $o^1$ and $o^2$.  For the boundary points $\gamma^1$ in $\pt_1$ and $\gamma^2$ in $\pt_2$, let $\gh_1=\gb_{\gamma^1}$ and $\gh_2=\gb_{\gamma^2}$ be the height functions associated with $(\T_1,o^1,\gamma^1)$ and $(\T_2,o^2,\gamma^2)$. The horospheric products $\T_1$ and $\T_2$, denoted by $\dl$ with the vertex set given by 
$$
\Big\{(x^1,x^2) \in \T_1 \times \T_2 \ : \  \gh_1(x^1)+\gh_2(x^2)=0 \Big\}.
$$
We declare vertices $ (x^1,y^1)$ and $(x^2,y^2)$ to be adjacent if $x^1$ and $y^1$ are adjacent in 
$\T_1$ and $x^2$ and $y^2$ are adjacent in $\T_2$.

An explicit description of the graph metric for regular trees is due to  Bertacchi \cite{Bertacchi}. An extension of her result that applies to all trees with no restriction is given by Kaimanovich and Sobiescky in \cite{Kaimanovich-Sobieczky}. In order to state the theorem, let $d, d_1,$ and $d_2$ denote the graph metrics for $\dl$, $\T_1$ and $\T_2$,
respectively. The distance between vertices $\bx=(x^1,x^2)$ and $\by=(y^1,y^2)$  is given by 
$\dl$, we have
$$
d(\bx,\by) = d_1(x^1,y^1)+d_2(x^2,y^2) - |\gh_1(x^1)-\gh_1(y^1)|.
$$

The following lemma is an immediate implication of this formula. 

\begin{lemma}\label{lem:dl-metric}
Keeping the notation as in the previous paragraph, then for all vertices $\bx=(x^1, x^2)$ and $\by= (y^1, y^2)$ of $\dl$ we have
\begin{equation}
\gb_{\bx}(\by)= \gb_{x^1}(y^1)+ \gb_{x^2}(y^2) - | \gh_1(x^1)- \gh_1(y^1)|+ | \gh_1(x^1)|.  
\end{equation} 

\end{lemma}

\subsection{\bf \bf The height compactifications of $\dl$} \label{sec:map-tree}
In this section we will introduce and study a particular mapping compactification of $\dl$, associated to the height function on $\dl$. 
Consider $\dl$ as a subgraph of the graph $\T_1 \times \T_2$, and hence as a subset 
of $ \ot_1 \times \ot_2$. Let $\overline{\ZZ}$ denote the two-point compactification of $\ZZ$ and consider the map
\[ h: \dl \to \overline{\ZZ} \]
defined by  $h(x^1,x^2)= \gh_1(x^1)$. Hence we obtain the map
$\Psi: \dl \to \overline{\T_1}  \times \overline{\T_2}\times \overline{\ZZ}$ defined by 
$$
\Psi(\bx) = (x^1,x^2, \gh_1(x^1) ).
$$

The {\it height compactification} of $\dl$ is defined 
to be the compactification of $\dl$ with respect to the map $\Psi$ and will be denoted by $\gcom$.

Our first goal is to give a precise characterization of $\gcom$ and its topology. 
 Let $\bx_n \in \dl$ be a convergent sequence in $\gcom$ and denote its limit by $(\bxi, \eta)$, where $\bxi \in
 \overline{\T_1} \times \overline{\T}_2 $ and $ \eta \in \overline{\ZZ}$. Writing $\bx_n=(x_n^1, x_n^2)$ and $\bxi= (\xi^1, \xi^2)$ this is equivalent to the convergence of $ x_n^i$ to $\xi^i$ in $\ot_i$ for $i=1,2$ and the convergence of
 $\gh_1(x_n^1)$ to $\eta$ in $ \overline{\ZZ}$.

For our characterization, we need to consider the following subsets of $\ot_1\times \ot_2 \times \ozz$:
\begin{equation}\label{decom:geometric}
\begin{split}      
\GC_1&=\pt_1\times\{\gamma^2\} \times \{+\infty\}, \\
\GC_2 &=  \{\gamma^1\} \times \pt_2 \times \{- \infty\}, \\
\GZ &= \{\gamma^1\} \times \{\gamma^2\} \times \ZZ,  \\
\GT_1 &=  \Big\{(x^1,\gamma^2,\gh_1(x^1) ) : x^1 \in \T_1 \Big\}\\
\GT_2 &=  \Big\{(\gamma^1,x^2,-\gh_2(x^2) : x^2 \in \T_2 \Big \}.
\end{split}
\end{equation}

\begin{theorem}\label{thm:geo} Assume that the degree of every vertex of $\T_1$ and $\T_2$ is at least $3$. With sets defined in \ref{decom:geometric}, we have 
\begin{enumerate}
\item The boundary of the compactification $\gcom$ is the union of $\GC_1$, $\GC_2$,
$\GT_1, \GT_2$, and $\GZ$.
\item $\GC_1$ and $\GC_2$ are compact
\item The set of limit points of $\GT_1$ and $\GT_2$ are equal to 
$\GC_1 \sqcup \GZ$ and $\GC_2 \sqcup \GZ$, respectively. 
\item The set of limit points of $\GZ$ consist of points $(\gamma^1,\gamma^2,
-\infty)$  and $(\gamma^1,\gamma^2,+\infty)$.
\end{enumerate}

\end{theorem}

\begin{proof}  
We will start by a brief description of convergent sequences in $\gcom$. Verifying the claims in each case is straightforward.
\begin{enumerate}
 \item (\emph{Eventually constant}.) For $\bx_n$ to converge to a point in $\dl$ it has to be eventually constant. The converse is clearly true. The limit will be the eventual value of $ (\bx_n, \gh_1(x_n^1))$. 
\item (\emph{Radial in the direction of $\gamma^2$}.) The sequence $\bx_n$ converges to a point of the form $$(\xi^1,\gamma^2,+\infty) \in \GC_1.$$ This happens if and only if $x_n^1 \to \xi^1 \in \overline{\T_1}$, $x_n^2 \to \gamma^2 \in \ot_2$, and $ \gh_1(x_n^1) \to + \infty$. 
\item (\emph{Radial in the direction of $\gamma^1$.})  The sequence  $x_n$ converges to a point of the form $$(\gamma^1,\xi^2,-\infty) \in \GC_2.$$
 This happens if and only if $x_n^1 \to \gamma^1 \in \ot_1$, $x_n^2 \to \xi^2 \in \ot_2$, 
 and $ \gh_2(x_n^2) \to + \infty$.
\item \emph{(Horocyclic.)} The sequence $\bx_n$ converges to 
$$
(\gamma^1,\gamma^2,k) \in \GZ.
$$
This is the case when 
$x_n^1 \to \gamma^1$ and  $x_n^2 \to \gamma^2$, while  $\gh_1(x^1_n) \to k$.  
\item \emph{(Going to infinity in the first component.)} The sequence $x_n$ converges to  
$$(\gamma^1,x^2,-\gh_2(x^2)) \in \GT_1.$$ For this to happen one requires $x_n^1 \to \gamma^1$ and $x^2 $ to be the eventual value of the sequence $x_n^2$. Note that this forces the sequence $-\gh_2(x_n^2)$ to be eventually constant.

\item (Going to infinity in the second component) The sequence $\bx_n$ converges to $$(x^1,\gamma^2,\gh_1(x^1)) \in \GT_1.$$ Similar to the previous case, this happens when  $x_n^2 \to \gamma^2$ in $\ot_2$ and $x^1 $ is  the eventual value of the sequence $x_n^1$. 
\end{enumerate}

One can easily verify that all these convergences are possible. For instance, 
for $x^2 \in \T_2$, using the fact that the degree of each vertex in both trees are at least three, one can find a divergent sequence of points $x_n^1 \in \T_1$ with $\gh_1(x_n^1)= -\gh_2(x^2)$. By Lemma~\ref{lem:bounded-tree}, we have  $x_n^1\to\gamma^1$ in $\ot_1$, hence $(x_n^1,x^2,\gh_1(x_n^1)) \to (\gamma^1,x^2, -\gh_2(x^2) )$ in $\gcom$.  

We shall now show that points in $\GT_1$ and $\GT_2$ can also be realized as limits of sequences in $\gcom$.  
Fix $\xi^2$ in $\pt_2 \setminus \{ \gamma^2 \} $ and pick a sequence  $(x^2_n)$ in $\T_2$ converging to $\xi^2$. Pick also a divergent sequence $(x_n^1) \in \T_1$ with the property that $-\gh_2(x_n^2)=\gh_1(x_n^1)$. Since, $x^2_n \to \xi^2$, we have $\gh_2(x^2_n) \to +\infty$,  implying that $\gh_1(x_n^1) \to -\infty$.  Hence,  $x_n^1 \to \gamma^1$. Thus, the sequence $(x_n^1,x_n^2)$ converges to $(\gamma^1,\xi^2,-\infty)$ in the height compactification of $\dl$. Now, by  choosing a sequence of boundary points $\xi^2_n$ that converges to $\gamma^2$ and a diagonal argument, one can construct a sequence of points converging to $(\gamma_1, \gamma_2,  - \infty)$. The case (3) can be dealt with similarly.

Finally we will show that $\GZ \subseteq \gcom$. 
Since every vertex in $\T_1$ and $\T_2$ has degree at least three, one can find for every $k \in \ZZ$ a divergent sequences  $\bx_n=(x_n^1,x_n^2) \in \dl$  such that  $\gh_1(x_n^1)=- \gh_2(x_n^2)=k$.  By Lemma~\ref{lem:bounded-tree}, $x_n^1\to\gamma^1$  and $x_n^2\to \gamma^2$, implying that $\bx_n$ converges to $(\gamma^1,\gamma^2,k)$ in $\gcom$. It is easy to see that as $k \to \pm \infty$  the sequence of points $(\gamma^1,\gamma^2,k)$ converge to $(\gamma^1,\gamma^2,+\infty)$ or  $(\gamma^1,\gamma^2,-\infty)$.

  \end{proof}
\section{\bf The Busemann compactification of $\dl$} \label{sec:Busemann}
In this section, we will {give a characterization of} the Busemann compactification of $\dl$. 
Along the way, we will show that a sequence converges in the Busemann compactification of $\dl$ if and only if it converges in the height compactification. This leads to an identification of $\gcom$ and $\bcom$. We will start by defining a family of real-valued functions on $\dl$. A point in $\dl$ is always denoted by $\bx=(x^1, x^2)$.

\begin{enumerate}
\item  For $\xi^1 \in \pt_1$, set $\gB_{\xi^1}(\bx)= \gb_{\xi^1}(x^1)$.
\item For $\xi^2 \in \pt_2$, set $\gB_{\xi^2}(\bx)= \gb_{\xi^2}(x^2)$. 
\item For $y^1 \in \T_1$, set $\gB_{y^1}(\bx)=\gb_{y^1}(x^1)+\gh_{2}(x^2)+\gc_{\gh_1(y^1)} (x^1)$.
\item For $y^2 \in \T_2$, set $\gB_{y^2}(\bx)=\gh_{1}(x^1)+\gb_{y^2}(x^2)+\gc_{-\gh_2(y^2)}(x^2)$.
\item For $k \in \ZZ$, set $\gB_{k}(\bx)= \mathfrak{c}_k(x^1)$.
\end{enumerate}

We will show $\bcom$ precisely consists of the functions $\gB_{\ast}$ defined above. Before we proceed, let us set for $i=1,2$ the following notation:
\begin{equation}\label{parts}
\BC_i= \{ \gB_\xi\ : \ \xi \in \pt_i\},\  
\BT_i = \{ \gB_y \ :\ y \in \T_i \}, \ 
\BZ = \{ \gB_k\ :\ k \in \ZZ\}.
\end{equation}
 
\begin{lemma}\label{lem:compactA}
For the sets defined in \eqref{parts}, we have
\begin{enumerate}
\item $\BC_1$ and $\BC_2$ are closed. 
\item The set of limit points of $\BZ$ consists of the points
$\gh_{1}$ and  $\gh_{2}$.
\item The set of limit points of $\BT_i$ is equal to $\BC_i \cup \BZ$ for $i=1,2$.

\end{enumerate} 
\end{lemma}
\begin{proof}
 It is easy to verify that $\BC_1$ and $\BC_2$ are closed. In order to 
determine the limit points of $\BZ$, fix $(x^1,x^2)$ in $\dl$.  
It is easy to see that for sufficiently large values of $k$ 
we have  $$\gc_k(x^1)=|k| - |k-\gh_1(x^1)|=\gh_{1}(x^1).$$ Likewise, as $k \to - \infty$, the eventual value of 
$\gc_k(x^1)$ will be $\gh_{2}(x^2)$. We deduce that the only limit points of $\BZ$ are $\gb_{\gamma^1}=\gh_1$ and $\gb_{\gamma^2}=\gh_2$. It remains to prove (3).  We will consider three cases. First, assume that $(x_n^1) $ is a sequence in $\T_1$ such that $x_n^1 \to \xi^1$ and $\gh_1(x_n^1) \to +\infty$.  This implies that $\gc_{\gh_1(x_n^1)}(x^1) \to \gh_1(x^1)$, which together with the
equality $\gh_2(x^2) = -\gh_1(x^1)$ implies that
$$
\lim_{n\to \infty}\gB_{x^1_n}(x^1,x^2) = \gb_{\xi^1}(x^1)+ \gh_{2}(x^2) +\gh_1(x^1)=  \gb_{\xi^1}(x^1).
$$
This, in particular, proves that the boundary point of $\BT_1$ includes $\BC_1$. Now
second, if $(x_n^1) $ is a sequence in $\T_1$ such that $\gh_1(x_n^1) \to -\infty$. By a similar argument, one can show that the only possible limit point is $\gh_1$. 
Finally, suppose that $(x_n^1) $ is a divergent sequence in $\T_1$ such that $\gh_1(x_n^1) $ is eventually constant $k$.  This implies that $x_n\to \gamma^1$ and $\gc_{\gh_1(x_n^1)}(x^1) \to \gc_k(x^1)$, which together with the
equality $\gh_2(x^2) = -\gh_1(x^1)$ implies that
$$
\lim_{n\to \infty}\gB_{x^1_n}(x^1,x^2) = \gh_1(x^1)+ \gh_{2}(x^2) +\gc_k(x^1)=  \gc_k(x^1)=\gB_k(x^1)\,
$$
which implies $\BZ$ is included in the boundary point of $\BT_1$. Therefore, the boundary point of $\BT_1$ is equal to  $\BC_1 \cup \BZ$.  The proof of the other case is similar. 
\end{proof}

The next theorem identifies the Busemann functions of $\dl$ that arise from divergent sequences of unbounded heights. 

\begin{theorem}[Limits of sequences of unbounded height]\label{thm:dl-non-constant}
Let $\bx_n=(x_n^1,x_n^2)$ be a divergent sequence of points of unbounded height in $\dl$ such that $\gb_{\bx_n}$ converges point-wise.  
Then the limiting function is either  of the form $\gB_{\xi^1}$ for some  $\xi^1 \in \pt_1 \setminus \{ \gamma^1 \}$ or $\gb_{\xi^2}$ for some  $\xi^2 \in \pt_2 \setminus \{ \gamma^2 \}$.
Moreover, $\bx_n$ converges in $\gcom$ in the former case to $(\gamma^1,\xi^2,-\infty)$ and
in the latter case to $(\xi^1,\gamma^2,+\infty)$.  
\end{theorem}
\begin{proof}

Let us first consider the case that $\gh_1(x_n^1) \to - \infty$. This implies that $x_n^1 \to \gamma^1$ and $\gb_{x_n^1}\to \gh_{1}$. 
By Lemma~\ref{lem:constant-cocycle}, for sufficiently large values of $n$ we have
$$
-|\gh_1(x_n^1)-\gh_1(x^1)|+ |\gh_1(x_n^1)| = -\gh_1(x^1).
$$
Applying Lemma~\ref{lem:dl-metric} to write the Busemann functions in the horospheric product in the terms of the one in the individual trees, we can write  for every $\by=(y^1,y^2)$ in $\dl$
$$
\gb_{\bx_n}(y^1,y^2) = \gb_{x^1_n}(y^1) + \gb_{x^2_n}(y^2) -|\gh_1(x_n^1)- \gh_1(y^1)|+ |\gh_1(x_n^1)|.
$$

Putting these together, it follows that $\gb_{x_n^2}$ also converges. Because the height of 
the sequence $(x_n^2)$ converges to $+ \infty$, $\gb_{x_n^2}$ must converge to a Busemann function 
$\gb_{\xi^2}$, where $\xi^2\not=\gamma^2$ is in the geometric boundary of $\T_2$. The limiting function will be
$
 \gb_{\xi^2} (y^2).
$

Observe that in this case, the convergence $\gb_{\bx_n} \to \gB_{\xi^2}$ in $\bcom$ corresponds to the convergence
$\Psi(\bx_n)$ to $(\gamma^1,\xi^2,-\infty)$ in $\gcom$. The case $\gh_1(x_n^1) \to + \infty$ can also be dealt with in a similar way.

\end{proof}

\begin{theorem}[Limits of sequences of bounded height]\label{thm:dl-constant}
Let $\bx_n$ be a divergent sequence of bounded height in $\dl$. Suppose, further, that the corresponding sequence of Busemann functions $\gb_{\bx_n}$ converges point-wise. Then the limiting function is exactly one of the following form:
\begin{enumerate}
\item $\gB_k$ for some $k \in \ZZ$. 
\item $\gB_{y^1}$ for some $y^1 \in \T_1$.
\item $\gB_{y^2}$ for some $y^2 \in \T_2$.
\end{enumerate}
Moreover, $\bx_n$ converges to $(y^1,\gamma^2,\gh_1(y^1))$ in (1), to $(\gamma^1,y^2,-\gh_2(y^2))$ in (2), and to $(\gamma^1,\gamma^2,k)$ in (3). 
\end{theorem}

\begin{proof}
Since $\bx_n$ diverges, so  $(x_n^1)$  or $(x_n^2)$  must be divergent. \\

{\it Case (1): divergence in both directions}.
First, assume that $(x_n^1)$  or $(x_n^2)$ both diverge. 
Since $|\gh_1(x_n^1)|=|\gh_2(x_n^1)|$ are assumed to be bounded by a constant $D$, it follows from Lemma~\ref{lem:bounded-tree} that $x_n^i\to \gamma^i$ for $i=1,2$, and $\gb_{x_n^i} $ converges point-wise for $i=1, 2$. By Lemma~\ref{lem:dl-metric} we have for $\by \in \dl$:
$$
\gb_{\bx_n}(y^1,y^2) = \gb_{x^1_n}(y^1) + \gb_{x^2_n}(y^2) -|\gh_1(x_n^1)- \gh_1(y^1)|+ |\gh_1(x_n^1)|.
$$
From the convergence of $\gb_{\bx_n}$, $\gb_{x_n^1}$, and $\gb_{x_n^2}$ we deduce that the sequence of functions
$$
\by \mapsto -|\gh_1(x_n^1)-\gh_1(y^1)|+ |\gh_1(x_n^1)|.
$$
also converges. Set $a_n = \gh_1(x_n^1)$. We will show that $a_n$ is eventually constant. Set $f_c(x)= |x|- |x-c|$. 
Using the fact that $|a_n | \le D$, we see that the sequence 
\[ f_{D} (a_n)+ f_{-D}(a_n)= 2|a_n|- | a_n - D| - |a_n+D|= 2|a_n |-2D \]
is eventually constant. Hence $a_n$ is either constant or takes two values $\pm k$, for some $k>0$. Note, however, that  $f_D(k)=k - (D-k)=2k-D$ while $f_D(-k)= -D$, implying that $k=0$, claiming that $a_n$ is eventually constant. 
Therefore,
 $\gh_1(x_n^1)$ is eventually constant $k$ and 
 $$
\lim_n\gb_{\bx_n}(y^1,y^2)  = \gb_{\gamma^1}(y^1) + \gb_{\gamma^2}(y^2)+
\gc_k(y^1).
 $$
Note that $\gb_{\gamma^1}(y^1) = \gh_1(y^1)$ and  $\gb_{\gamma^2}(x^2)=\gh_2(y^2)=-\gh_1(y^1)$, therefore,
$$
\lim_n\gb_{\bx_n}(y^1,y^2) = \gB_k(y^1)= \gc_k(y^1).
$$

{\it Case (2): divergence in $\T_1$ component}. Now assume the sequence $(x_n^1)$ is divergent in $\T_1$, but the sequence $(x_n^2)$ is not divergent in $\T_2$. Hence, the sequence $(x_n^2)$ is eventually equal to $x^2$. The rest of the argument is similar to Case (1),  
 $$
\lim_n\gb_{\bx_n}(y^1,y^2)   = \gb_{\gamma^1}(y^1) + \gb_{x^2}(y^2)+
\gc_k(x^1),
 $$
where $\gh_2(x^2)=-k$. \\

{\it Case (3): divergence in $\T_2$ component}

 This is when $(x_n^1)$ is eventually $x^1$ and $(x_n^2)$ is divergent, the proof is similar to previous cases, and 
  $$
\lim_n\gb_{\bx_n}(y^1,y^2)   =  \gb_{x^1}(y^1)+ \gb_{\gamma^2}(y^2) +
\gc_k(x^1),
 $$
where $\gh_2(x^1)=k$.

\end{proof}

Putting Lemma~\ref{lem:compactA}, Theorems~\ref{thm:dl-constant} and ~\ref{thm:dl-non-constant} together, we can conclude that the boundary of Busemann compactification of $\dl$ is exactly equal to $ \BC_1 \sqcup  \BC_2 \sqcup \BT_1 \sqcup \BT_2 \sqcup \BZ $.
\begin{theorem}
Let $\T_1$ and $\T_2$ be two trees without any degree two or one vertices. Then, the boundary points of the Busemann compactification $\bcom$ is equal to
$$
 \BC_1 \sqcup  \BC_2 \sqcup \BT_1 \sqcup \BT_2 \sqcup \BZ. 
 $$
  
\end{theorem}

\begin{theorem}\label{thm:bus=HM}
Let $\T_1$ and $\T_2$ be two trees without any degree two or one vertices. Then, the Busemann compactification of $\dl$   is isomorphic to the  height compactification of
$\dl$.
\end{theorem}
\begin{proof}
Define map $\Theta: \Psi(\dl) \to \Phi(\dl) $ as $\Theta( \Psi(\bx) ) = \mathfrak{b}_\bx$.  By Theorems~\ref{thm:dl-constant} and \ref{thm:dl-non-constant}, we have 
$\bx_n=(x_n^1,x_n^2) $ converges to $\gB$ in the Busemann compactification if and only if 
$(x_n^1,x_n^2,\gh_1(x_n^1)) $ converges in the height compactification.  We can extend $\Theta$ continuously  from $\bcom$ to $\gcom$.
\end{proof}

\subsection{Trees with degree $2$ vertices}
Let us now consider the case that  $\T_1$ or $\T_2$ includes vertices of degree $2$. 
We will show that the Busemann compactification and the height compactification of $\dl$ need not be isomorphic. For example, assume that $\T_2$ is an bi-infinite path (that is, isomorphic to the standard Cayley graph of $\ZZ$). It is easy to see that $\dl$ is isomorphic to $\T_1$. However, one can see that if the degree of every vertex of $\T_1$ is at least $3$, then the Busemann compactification and the height compactification
of $\T_1$ are not isomorphic. Indeed, corresponding to every $k \in \ZZ$, by choosing a divergent sequence of vertices in $\T_1$ with $\gh_1(x_n)=k$, one obtains a point in the height compactification, yielding an additional copy of $\ZZ$ in $ \gcom$.

In order to clarify the role of such sequences, let $(\T, o, \gamma)$ be a rooted tree pointed at the fixed geodesic $\gamma$. Assume, further, that the degree of every vertex is at least $2$. Denote by $F(\T)$ the set of all $k \in \ZZ$  such that 
$ H_k(\T)$ is infinite.
For $k \in F(\T)$,  there exists a ``horocyclic''  sequence of $(x_n)$ of level $k$ such that converges to $\gamma$ in the geometric boundary of $\T$.

Now, if $k \in F(\T_1)$, it follows from the proof of Theorem~\ref{thm:dl-constant} that the function
$$\gB_{y^1}(\bx)=\gb_{y^1}(x^1)+\gh_{2}(x^2)+\gc_{k} (x^1)$$
belongs to the Busemann compactification of $\dl$ if and only if $-k \in F(\T_2)$. 
More generally, the limit points of $\bcom$ are exactly functions of the following kind:
\begin{enumerate}
\item For $\xi^1 \in \pt_1$, set $\gB_{\xi^1}(\bx)= \gb_{\xi^1}(x^1)$.
\item For $\xi^2 \in \pt_2$, set $\gB_{\xi^2}(\bx)= \gb_{\xi^2}(x^2)$. 
\item For $y^1 \in \T_1$ and $-\gh_1(y^1) \in F(\T_2)$, set 
$$
\gB_{y^1}(\bx)=\gb_{y^1}(x^1)+\gh_{2}(x^2)+\gc_{\gh_1(y^1)} (x^1).
$$
\item For $y^2 \in \T_2$ and $-\gh_2(y^2) \in F(\T_1)$, set 
$$
\gB_{y^2}(\bx)=\gh_{1}(x^1)+\gb_{y^2}(x^2)+\gc_{-\gh_2(y^2)} (x^2).
$$
\item For $k \in F(\T_1)$ and $-k \in F(\T_2)$, set $\gB_{k}(\bx)= \mathfrak{c}_k(x^1)$.
\end{enumerate}

The proof of the following theorem is along the same lines as the proof of Theorems
~\ref{thm:dl-constant} and~\ref{thm:dl-non-constant}.

\begin{theorem}\label{thm:general}
Let $\T_1$ and $\T_2$ be two infinite trees whose vertices have degree at least two.  Then, the divergent  sequence $\bx_n=(x^1_n,x_n^2)$ converges in the Busemann compactification of $\dl$ if and only if $\bx_n$ converges in the height compactification of $\dl$ and $\lim_n \gh_1(x_n^1) \in F(\T_1)\cup\{+\infty,-\infty\}$ and and $\lim_n \gh_2(x_n^2) \in F(\T_2)\cup\{+\infty,-\infty\}$. 

\end{theorem}

\begin{remark}\label{cocycle}

The characterization of the Busemann compactification for $\dl$ in connection with 
the general law of large number proved in \cite{Karlsson-Ledrappier-law} provides interesting information about the asymptotic behavior of certain stochastic processes.  Before we describe this application, let us  recall the setting of this theorem. 
Let $(X,d)$ be a proper metric space with a base point $o$.  Denote  the group of isometries of $X$ by $\textrm{Iso}(X)$.
Let $(\Omega,\mu)$ be a standard probability space with an ergodic measure preserving transformation $L:\Omega \to \Omega$. Let $g: \Omega \to \textrm{Iso}(X)$ be a measurable map. We define the associated ergodic cocycle $Z_n$ by
$$
Z_n(\omega) = g(\omega) g(L\omega) \cdots g(L^{n-1}\omega).
$$
We write $Z_no:=Z_n(\omega) o$ when $\omega$ is fixed.

 Then \cite[Theorem 1.1]{Karlsson-Ledrappier-law} states that when $Z_n$ is integrable there exists a measurable map assigning 
to  almost every $ \omega \in \Omega$ a Busemann function $\gB$ (depending on $\omega$) such that 
\[ \lim_{n \to \infty} - \frac{1}{n} \gB ( Z_n o) = A, \]
where $ A:= \displaystyle\lim_{n \to \infty} \dfrac{1}{n} d(Z_n o , o)$  by Kingsman's subadditive ergodic theorem. In other words, the direction at which the process $Z_no$ diverges can be detected
by a Busemann function.

When the Busemann compactification of $X$ can be geometrically described, this theorem can provide 
explicit information about the typical asymptotic behavior of $ Z_n o$ as $n \to \infty$. Let us first consider 
the case that $X$ is a locally infinite tree with a base point $o$, and $Z_n$ be as above. 

We assume that $A>0$. In view of Karlsson-Ledrappier's law of large numbers \cite{Karlsson-Ledrappier-law} and the fact that the Busemann functions for a tree are given by the points on the geometric boundary, there exists a boundary point $\gamma$ such that 
$$
\lim_{n\to\infty} \frac{\gB_\gamma(Z_n o)}{n} = -A.
$$
We now claim that for every  $\gamma' \neq \gamma$ we have 
\[ \lim_{n\to\infty} \frac{\gB_{\gamma'}(Z_n o)}{n} = A. \]

First, note that since $ Z_n o $ converges to $\gamma$, hence for every $\gamma' \neq \gamma$, for sufficiently large values of $n$, 
the value of $ d(o,[o,Z_no]_{\gamma'}) $ will be a constant $C$ depending on $\gamma$ and $\gamma'$.  
This implies that 
 
\[ \gB_{\gamma'}(Z_n o) = d(o,Z_n o ) -2 d(o,[o,Z_n o]_{\gamma'}) = d(o,Z_n o ) - 2C.  \]
it follows that 
\begin{equation}\label{upper}
 \lim_{n\to\infty} \frac{\gB_{\gamma'}(Z_n o)}{n} =  \lim_{n\to\infty} \frac{d(Z_n o, o)}{n} = A.
\end{equation}
The claim follows  immediately.

Now, suppose $\T_1$ and $\T_2$ are trees of degree at least $3$. 
Let $Z_n=(X_n, Y_n)$ be an ergodic integrable cocycle taking values in the affine group of $\dl$. By this we mean that $X_n$ and $Y_n$ belong to isometry groups of 
$\T_1$ and $\T_2$, respectively, and that for all $n \ge 1$, one has $X_n \gamma_1=\gamma_1$,  $Y_n \gamma_2 = \gamma_2$, and 
\begin{equation}\label{sum}
\gh_1(X_no^1) + \gh_2(Y_n o^2) =0.
\end{equation}
This implies that $X_n$ and $Y_n$ are, respectively, integrable ergodic cocycles in the group of isometries of $ \T_1$ and $\T_2$. We will remark that when $\T_1$ and $\T_2$ are non-isomorphic homogenous trees, then the affine group of $\dl$ is a non-discrete locally compact group and coincides with 
the full isometry group of $\dl$. In the case that $\T_1$ and $\T_2$ are isomorphic homogenous trees, then the affine group of $\dl$ has index $2$ in the 
full isometry group of $\dl$. In either case, one can easily construct many cocycles (e.g. coming from random walks) of the form $Z_n=(X_n, Y_n)$ of the above form. 

By the distance formula in $\dl$ from \cite{Kaimanovich-Sobieczky} we have 
$$
d(Z_no,o) = d(X_no^1,o^1) + d(Y_n o^2,o^2) - | \gh_1(X_n o^1) |.
$$
Suppose $ \displaystyle \lim_{n \to \infty} \dfrac{1}{n} d(Z_n o , o)=A>0$. Then at least one of $\displaystyle \lim_{n \to \infty} \dfrac{1}{n} d(X_n o^1 , o^1)$ or $\displaystyle \lim_{n \to \infty} \dfrac{1}{n} d(Y_n o^2 , o^2)$ must be positive. Without loss of generality, suppose that 
$\displaystyle \lim_{n \to \infty} \dfrac{1}{n} d(X_n o^1 , o^1)>0$.  We claim that exactly one of 
$$\gh_1(X_n o^1) \to + \infty \quad \textrm{ and }  \quad \gh_2(Y_n o^1) \to + \infty$$ holds.  Indeed, in view of the above discussion for the case of trees, one has that if 
$\gh_1(X_n o^1) \to + \infty$ does not hold, then $\gh_1(X_n o^1) \to -\infty$. Now, \eqref{sum} implies that  $\gh_2(Y_n o^1) \to + \infty$.
Suppose that $\gh_1(X_n o^1) \to + \infty$. Then $ d(X_no^1,o^1)  - | \gh_1(X_n o^1) |$ remains bounded as $n \to \infty$ and hence
$$ A:= \lim_{n\to\infty} \frac{d(Z_no,o)}{n} =  \lim_{n\to\infty} \frac{d(Y_no^2,o^2)}{n} = 
\lim_{n\to\infty} \frac{-\gh_2(Y_no^2)}{n}.
$$
A similar results holds when $\gh_2(Y_n o^2) \to + \infty$.
As a result, we see that in each case the Busemann function  appearing in Karlsson-Ledrappier's law of large numbers can be taken to be $\gh_1$ or $\gh_2$.

\end{remark}

\subsection*{Data Availability}
No datasets were generated or analyzed during the current study.
\bibliographystyle{amsalpha}
\bibliography{bib-PB}
\end{document}